\documentclass[twocolumn,         
               showpacs,longbibliography,            
               showkeys,preprintnumbers,     
               aps,                 
               prd,          	    
               letterpaper,             
               superscriptaddress,      
               nofootinbib,         
               tightenlines,        
               floats,floatfix      
               ]{revtex4-1}
\usepackage{physics}
\usepackage{epsf}
\usepackage{graphicx,color}
\usepackage{subfigure}
\usepackage{latexsym}
\usepackage{amsmath,amssymb}        
\usepackage[colorlinks=true,linkcolor=blue,citecolor=blue]{hyperref}
\usepackage{mathrsfs}
\usepackage{comment}
\usepackage{soul}
\usepackage{cancel}
\usepackage{framed}
\definecolor{purple}{rgb}{0.58,0.0,0.83}
\definecolor{orange}{rgb}{1,0.5,0}
\DeclareSymbolFontAlphabet{\mathrsfs}{rsfs}
\DeclareMathAlphabet{\mathcal}{OMS}{cmsy}{m}{n}

\begin{document}


\title{Numerical Solution Partial Differential Equations using the Discrete Fourier Transform}


\author{Daniela Estefan\'ia Rodr\'iguez-Lara}
\email{1130051c@umich.mx}
\affiliation{Facultad de Ciencias F\'isico Matem\'aticas, Universidad
              Michoacana de San Nicol\'as de Hidalgo. Edificio ALFA, Cd.
              Universitaria, 58040 Morelia, Michoac\'{a}n,
              M\'{e}xico.}

\author{Iv\'an  \'Alvarez-Rios}
\email{ivan.alvarez@umich.mx}
\affiliation{Instituto de F\'{\i}sica y Matem\'{a}ticas, Universidad
              Michoacana de San Nicol\'as de Hidalgo. Edificio C-3, Cd.
              Universitaria, 58040 Morelia, Michoac\'{a}n,
              M\'{e}xico.}

\author{Francisco S Guzm\'an}
\email{francisco.s.guzman@umich.mx}
\affiliation{Instituto de F\'{\i}sica y Matem\'{a}ticas, Universidad
              Michoacana de San Nicol\'as de Hidalgo. Edificio C-3, Cd.
              Universitaria, 58040 Morelia, Michoac\'{a}n,
              M\'{e}xico.}


\date{\today}


\begin{abstract}
In this paper we explain how to use the Fast Fourier Transform (FFT) to solve partial differential equations (PDEs). We start by defining appropriate discrete domains in  coordinate and frequency domains. Then describe the main limitation of the method arising from the Sampling Theorem, which defines the critical Nyquist frequency  and the aliasing effect. We then define the Fourier Transform (FT) and the FFT  in a way that can be implemented in one and more dimensions. Finally, we show how to apply the FFT in the solution of PDEs related to problems involving two spatial dimensions, specifically the Poisson equation, the diffusion equation and the wave equation for elliptic, parabolic and hyperbolic cases respectively.
\end{abstract}


\pacs{keywords: Partial Differential Equations of Physics -- Numerical Methods -- Fourier Transform}


\maketitle

\section{Introduction}
\label{sec:intro}

The behavior of physical systems is mostly described with Partial Differential Equations (PDEs). In some simple cases, when symmetries are assumed, closed solutions to such PDEs can be found, however in most cases, where no symmetries can be assumed, it is necessary to calculate numerical solutions. This is why the numerical solution of PDEs can be a discipline in itself, helpful at solving Physics problem of current interest. Inherent to numerical solution of PDEs is a variety of numerical methods used, among which there are spectral methods, methods based on finite differences, methods specific for elliptic problems or methods designed for initial value problems involving evolution (see e.g. \cite{guzman2023}).

In this journal various numerical methods to solve PDEs have been explained, for example for the Schr\"odinger equation \cite{RMFSchro} and wave equation \cite{guzman2010} using Finite Differences methods along with different evolution schemes, General Relativity Equations using Spectral Methods \cite{PseudoGR}, Relativistic Hydrodynamics Equations using Finite Volume based methods \cite{SphericalHydro,ShockTube}, including the evolution of the accretion of a fluid into a black hole \cite{SphericalAccretion}.

The aim of this paper is to explain the use of the Discrete Fourier Transform (DFT) and its use in the solution of PDEs, which would cover another useful numerical method. In the continuum, the method is easy to describe, and consists in transforming a PDE from the space-time domain into the frequency domain using the FT, which simplifies the problem by converting derivatives into algebraic terms, which are easier to handle and allowing an easy-to-obtain  solution to the PDF in the frequency domain. The solution to the problem in the space-time domain would simply be the inverse Fourier Transform of the solution found in the frequency domain. 

However, in practice, for rather general and complex to solve PDEs, a numerical solution is searched on a discrete domain, both, in the space and frequency domains and therefore the discrete version of the FT has to be defined, as well as its inverse. The discretization implies limitations on the spectrum of functions in the frequency domain due to sampling errors that have to be taken into account. In this paper we explain how to calculate the DFT, its inverse, as well as their limitations, followed by the description of the most popular method to calculate the FT, the Fast Fourier Transform (FFT) method  \cite{cooley1965algorithm} that speeds up the calculation.
We then proceed to explain, with examples of PDEs defined on domains with two spatial dimensions, the use of the FFT to solve elliptic, parabolic and hyperbolic equations, that we illustrate with the Poisson, Dissipation and Wave equations respectively.

The paper is organized as follows, in Section \ref{sec:method} we describe the DFT and its inverse, later in Section \ref{sec:applications} we show how to solve PDEs and finally in Section \ref{sec:conclusions} we draw some final comments.


\section{The Method}
\label{sec:method}

In this section, we describe the approximation of the Fourier Transform, starting with the one-dimensional case and extending the method to higher dimensions. These approximations rely on defining functions over a discrete domain and finding arithmetic representations of the Fourier operators on that domain. Before delving into this, however, it is crucial to introduce a key theorem in signal theory: the sampling theorem.

\subsection{Sampling Theorem and Its Implications}

Consider a continuous function \( p(t) \) sampled at regular intervals of length \( \Delta t \). If its FT \cite{spiegel1974schaum,press1986numerical}

\begin{equation}
    \mathcal{F}\left\lbrace p(t)\right\rbrace :=\hat{p}(f) = \int_{-\infty}^{\infty} p(t) e^{2\pi \mathrm{i} f t} \, dt,
    \label{eq:FT}
\end{equation}

\noindent has a finite bandwidth, meaning that \( \hat{p}(f) = 0 \) for frequencies higher than a critical one \( |f| \geq f_c \), where \( f_c = \dfrac{1}{2 \Delta t} \) is the Nyquist frequency, then the function \( p(t) \) is completely determined by its sampled values \( p_n = p(t_n) \) with \( t_n = n \Delta t \), for \( n = -\infty, \dots, -1, 0, 1, \dots, \infty \). In this case, \( p(t) \) can be reconstructed exactly from its samples using the following interpolation formula:

\begin{equation}
    p(t) = \Delta t \sum_{n=-\infty}^{\infty} p_n \dfrac{\sin\left[ 2\pi f_c (t - n\Delta t) \right]}{\pi (t - n\Delta t)}.
    \label{eq:p_Nyquist}
\end{equation}

\noindent The proof of this result involves the FT of both the Dirac comb and the rectangle function \citep{press1986numerical}. While we do not focus on the detailed proof here, we will discuss the key implications of this theorem for sampled signals.

First, any function that is constructed from discrete samples taken at regular intervals \( \Delta t \) will have a Fourier spectrum limited to the bandwidth \( |f| < f_c \). Second, if the sampling interval \( \Delta t \) is too large, a phenomenon known as aliasing will occur. Aliasing causes frequencies higher than the Nyquist frequency to be misrepresented, or  ``folded'', into the interval \( |f| < f_c \), leading to the distortion of the sampled signal.

With this in mind, any numerical or observational method that relies on a finite sample of a function can only capture a limited spectrum of frequencies. As a consequence, some frequencies of interest outside this spectrum may not be detected. One strategy to mitigate aliasing is to use multiple sets of samples with smaller intervals, ensuring that the signals are sampled consistently and that higher frequencies are accurately represented.

With this understood we are now prepared to apply a discrete version of the FT (\ref{eq:FT}).

\subsection{Discrete Fourier Transform}

To construct the discrete version of the FT, we first define a discrete, finite time domain:

\begin{equation}
    t_n = n\Delta t, \qquad n = 0,1,2,...,N-1,
    \label{eq:discrete_temporal_domain}
\end{equation}

\noindent where \(N\) is the number of points in the domain, and a function \(p(t)\) evaluated at these points is denoted as \(p_n\). With this, the integral in equation (\ref{eq:FT}) can be approximated by a Riemann sum:

\begin{equation}
    \hat{p}(f) = \sum_{j=0}^{N-1} p_j e^{2\pi \mathrm{i} f j \Delta t} \Delta t + \mathcal{O}(\Delta t^2). 
\end{equation}

\noindent Next, to establish a correspondence between the time and frequency domains, we define a discrete frequency domain instead of treating the frequency \(f\) as continuous:

\begin{equation}
    f_k = \dfrac{k}{N\Delta t}, \qquad k = 0,1,2,...,N-1,
    \label{eq:discrete_fourier_domain}
\end{equation}

\noindent where the FT is defined only at the points of this discrete frequency domain and denoted as \(\hat{p}_k = \hat{p}(f_k)\). Consequently, the FT takes the form

\begin{equation}
    \hat{p}_k = \sum_{j=0}^{N-1} p_j e^{2\pi \mathrm{i} j k / N} \Delta t + \mathcal{O}(\Delta t^2). 
\end{equation}

\noindent In this way, the DFT of the discrete function \(\Vec{p} = (p_0, p_1, p_2, ..., p_{N-1})^T\) is defined as

\begin{equation}
P_k = \sum_{j=0}^{N-1} p_j \omega_N^{jk}, 
\label{eq:DFT}
\end{equation}

\noindent where \(\omega_N := e^{2\pi \mathrm{i} / N}\). Equivalently, we can express this formula in matrix form:

\begin{equation}
    DFT(\Vec{p}) := \Vec{P} = \mathbf{W} \Vec{p},
    \label{eq:DFT2}
\end{equation}

\noindent where the matrix \(\mathbf{W}\) is given by

\begin{equation}
    \mathbf{W} = \begin{pmatrix}
    \omega_N^{0\cdot 0} & \omega_N^{0\cdot 1} & \cdots & \omega_N^{0\cdot (N-1)} \\
    \omega_N^{1\cdot 0} & \omega_N^{1\cdot 1} & \cdots & \omega_N^{1\cdot (N-1)} \\
    \vdots & \vdots & \ddots & \vdots  \\
    \omega_N^{(N-1)\cdot 0} & \omega_N^{(N-1)\cdot 1} & \cdots & \omega_N^{(N-1)\cdot (N-1)} \\
    \end{pmatrix}.
\end{equation}

\noindent where the dot in the exponents is a product. From these two expressions of the DFT, two important properties of the DFT arise:

\begin{enumerate}
    \item The DFT is periodic, i.e., \(P_{k+N} = P_k\).
    \item The inverse of the matrix \(\mathbf{W}\) is given by \(\mathbf{W}^{-1} = \dfrac{1}{N} \mathbf{W}^\dagger\), where \(\mathbf{W}^\dagger\) is the conjugate transpose of \(\mathbf{W}\).
\end{enumerate}

\noindent The second property allows us to recover the original function from its Fourier components via the inverse DFT ($iDFT$) given by:

\begin{equation}
    p_j = \dfrac{1}{N} \sum_{k=0}^{N-1} P_k \omega_N^{-jk},
\end{equation}

\noindent or equivalently, in matrix form:

\begin{equation}
    iDFT(\Vec{P}) := \Vec{p} = \dfrac{1}{N} DFT(\Vec{P}^*)^* = \dfrac{1}{N} \mathbf{W}^* \Vec{P},
\end{equation}

\noindent where ${}^*$ denotes the complex conjugate. As an example of the approximated FT using the DFT, consider the function

\begin{equation}
    p(t) = e^{-t^2},
\end{equation}

\noindent whose exact FT according to (\ref{eq:FT}) is

\begin{equation}
    \hat{p}(f) = \sqrt{\pi} e^{-(\pi f)^2}.
\end{equation}

\noindent We compute the DFT approximation in the discrete domain \([0,20)\), with a time step \(\Delta t = 20/N\), where \(N = 32\) points.

\begin{figure}
    \centering
    \includegraphics[width=8cm]{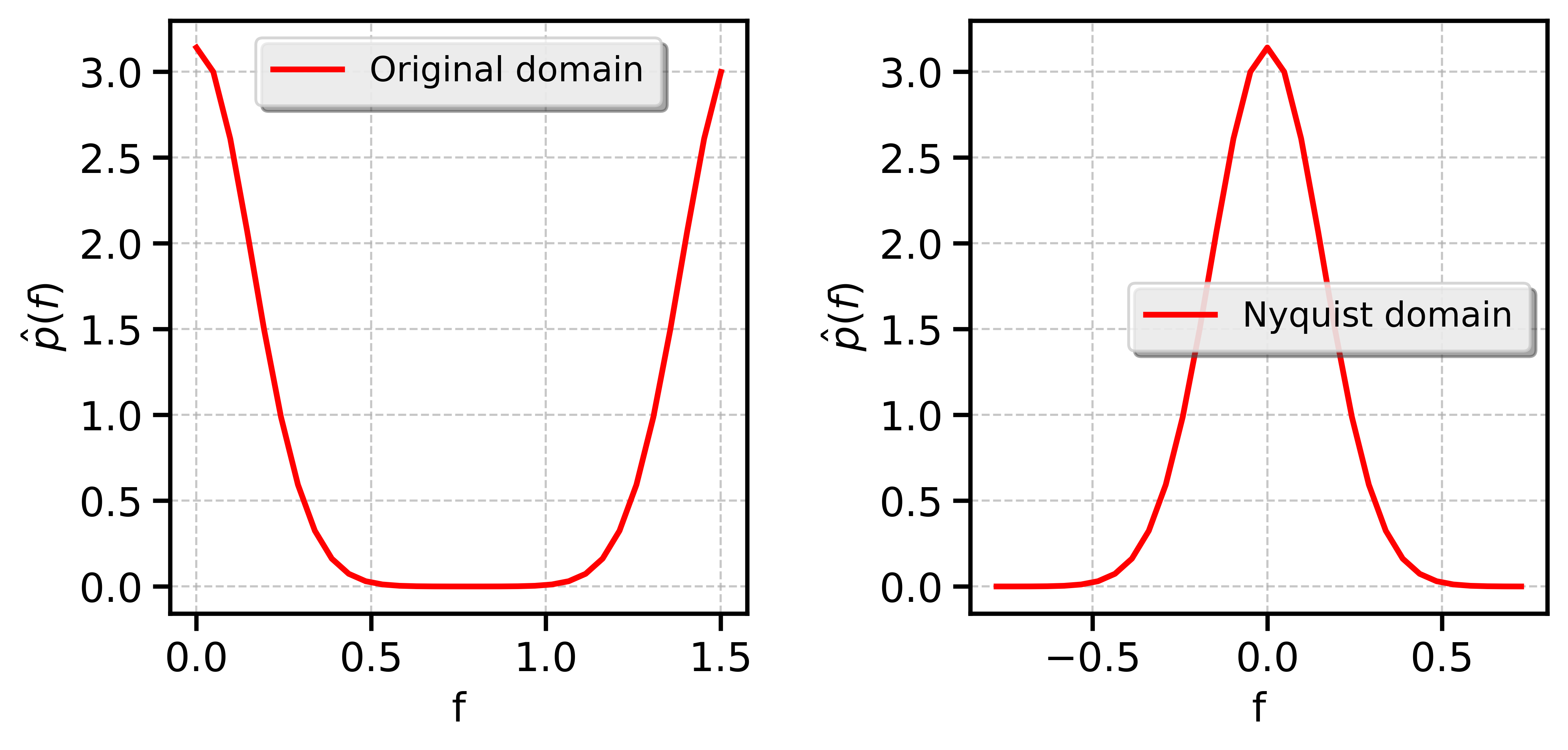}
    \caption{Left: Approximation of the FT in the original frequency domain. Right: The same signal after frequency shifting, aligned with the Nyquist frequency.}
    \label{f:FreqShift}
\end{figure}

An important aspect of spectral analysis is frequency shifting. The frequencies defined by the DFT in equation (\ref{eq:DFT}) do not correspond directly to the Nyquist frequency for indices \(k > N/2\). However, taking advantage of the signal's periodicity, we can shift these frequencies to align with the Nyquist frequency. This procedure is illustrated in Figure \ref{f:FreqShift}, where the left panel presents the Fourier transform in the frequency domain as obtained from formula (\ref{eq:DFT}), and the right panel shows the same signal shifted to account for periodicity, this time correctly aligned with the Nyquist frequency.

It is important to note that calculating the DFT requires \(\mathcal{O}(N^2)\) operations, making it computationally expensive. In order to optimize computer resources, Cooley and Tukey developed a more efficient algorithm, nowadays known as the Fast Fourier Transform (FFT) \cite{cooley1965algorithm}, which significantly reduces the computational cost and that we describe now.

\subsection{Fast Fourier Transform}

The essential idea behind the Fast Fourier Transform (FFT) is that the Discrete Fourier Transform (DFT) for a dataset with \(N\) points can be efficiently decomposed into two smaller DFTs, each with \(N/2\) points. One of these smaller DFTs involves the even-indexed points, while the other involves the odd-indexed points.

Given a discrete sequence \(p_j\), where \(j = 0,1,2,...,N-1\) and assuming \(N\) is even, the Fourier transform \(P_k\) can be expressed as the sum of two transforms, each of length \(N/2\):

\begin{eqnarray}
 P_k &=&  \sum_{j=0}^{N/2-1}p_{2j}\omega_N^{(2j)k} + \sum_{j=0}^{N/2-1} p_{2j+1}\omega_N^{(2j+1)k}\nonumber\\
 &=& \sum_{j=0}^{N/2-1}p_{2j}\omega_{N/2}^{jk} + \sum_{j=0}^{N/2-1} p_{2j+1}\omega_{N/2}^{jk} \omega_N^k.\nonumber
\end{eqnarray}

\noindent Here, \(\omega_N^k\) is independent of \(j\), so it can be factored out of the second summation. Therefore, we can rewrite the equation as:

\begin{equation}
P_k = P_k^E + \omega_N^k P_k^O, \quad k=0,1,2,\dots,N/2-1,
\label{eq:FFTpairodd}
\end{equation}

\noindent where \(P_k^E\) and \(P_k^O\) represent the Fourier transforms of the even and odd-indexed components of the sequence, respectively. In this decomposition, \(k\) runs over \(N/2\) values for both the even and odd transforms, effectively reducing the original problem into two smaller DFTs.

For these smaller DFTs, the following periodicity conditions apply:

\begin{eqnarray}
P^{E}_{k+N/2}&=&P^{E}_k,\nonumber\\
P^{O}_{k+N/2}&=&P^{O}_k.\nonumber
\end{eqnarray}

\noindent Additionally, the factor \(\omega^k_N\) satisfies the important identity:

\begin{equation}
    \omega^{k+N/2}_N = -\omega^k_N.
\end{equation}

\noindent The remaining \(N/2\) terms can be computed using:

\begin{equation}
P_{k+N/2} = P_k^E - \omega_N^k P_k^O, \quad k=0,1,2,\dots,N/2-1.
\label{eq:FFTpairodd_2}
\end{equation}

\noindent This process effectively reduces the computational complexity from \(\mathcal{O}(N^2)\) in the case of a direct DFT to \(\mathcal{O}(N \log N)\), making the FFT much more efficient for large \(N\).

\begin{figure}
    \centering
    \includegraphics[width=8cm]{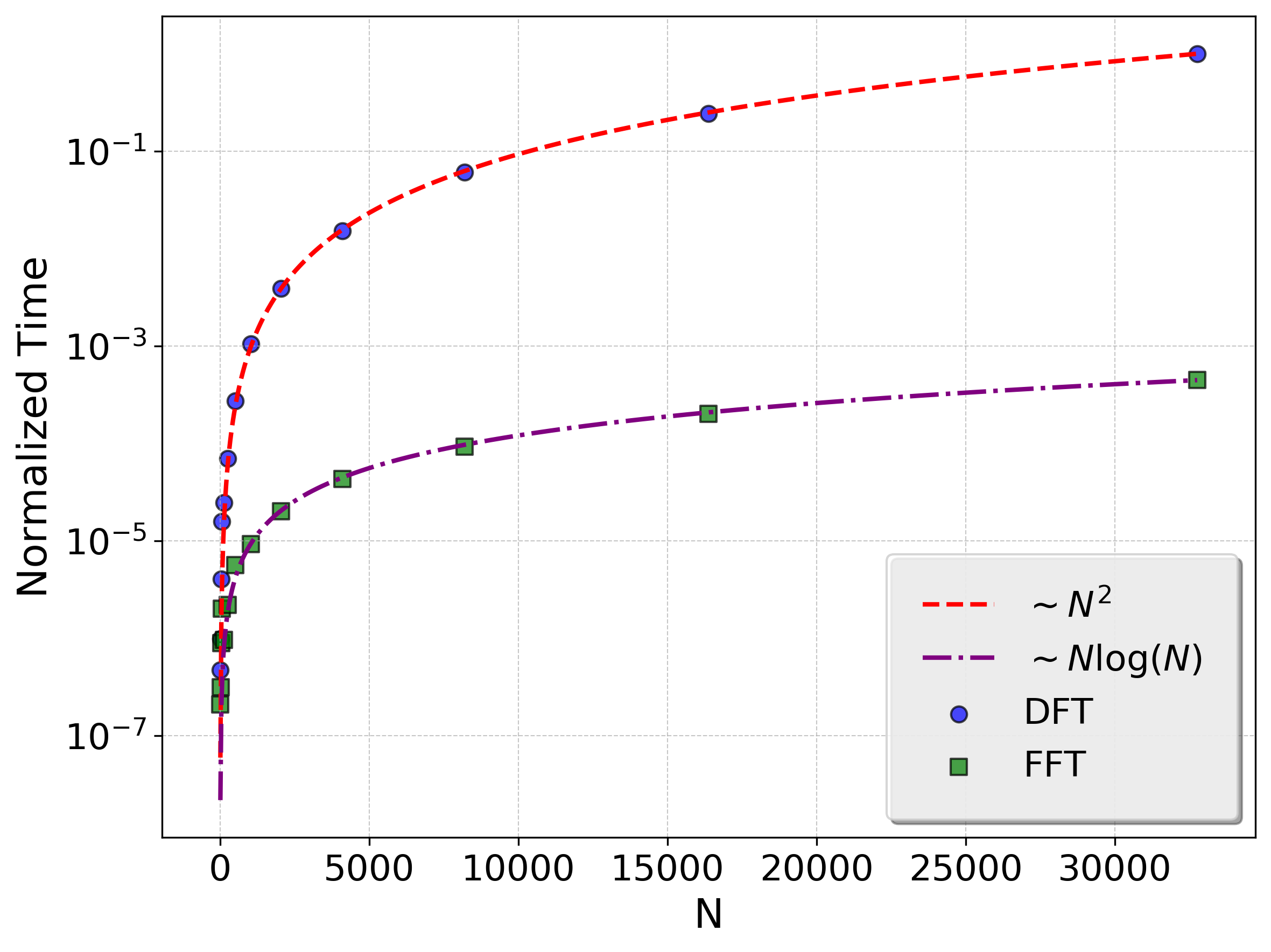}
    \caption{Normalized CPU time as a function of \(N\) for the calculation of the Discrete Fourier Transform using the DFT and FFT methods. Fits are also shown for the \(\sim N^2\) and \(\sim N \log N\) scaling.}
    \label{f:CPUtime}
\end{figure}

To illustrate this efficiency gain, we compare the performance of the DFT and the FFT for the previous example, using a time step \(\Delta t = 20 / N\), where \(N = 2^m\) and \(m = 3,4,5,\dots,15\). Figure \ref{f:CPUtime} shows the normalized CPU time as a function of \(N\), relative to the maximum value calculated for \(N=2^{15}\). As expected, the DFT scales as \(\sim N^2\), while the FFT scales as \(\sim N \log N\). For instance, when \(N = 2^{15}\), the FFT is approximately 2000 times faster than the DFT.

The number of operations is the main reason why the FFT is used in practice, and the fact that $N$ has to be of the form $2^m$ is a good price to pay. In the examples developed below we use the FFT approach.

\subsection{Discrete Fourier Transform in Two and More Dimensions}

Consider the discrete two-dimensional domain $D_d=\{(x_j,y_k) \}$ such that
\(x_j = j \Delta x\) and \(y_k = k \Delta y\) for \(j = 0, 1, \ldots, N_x-1\) y \(k = 0, 1, \ldots, N_y-1\). There we define the complex function \(p_{jk} = p(x_j, y_k)\) at each point of $D_d$, which in general can be complex. The two-dimensional Discrete Fourier Transform (DFT) of \(p_{jk}\), denoted as \(P_{ab}\), can be expressed as:

\begin{equation}
 P_{ab} = \sum_{j=0}^{N_x-1} \sum_{k=0}^{N_y-1} p_{jk} e^{2\pi \mathrm{i} aj/N_x} e^{2\pi \mathrm{i} bk/N_y}.
 \label{2D_DFT}
\end{equation}

\noindent The corresponding discrete frequency domains are defined as follows:

\begin{equation}
 \begin{matrix}
 k_a^x = \dfrac{a}{N_x \Delta x}, \quad a = 0, 1, 2, \ldots, N_x-1, \\
 k_b^y = \dfrac{b}{N_y \Delta y}, \quad b = 0, 1, 2, \ldots, N_y-1.
 \end{matrix}
\end{equation}

\noindent along each direction in the two dimensional frequency domain. In order to compute the DFT we can reorganize equation (\ref{2D_DFT}) in the following form:

\begin{equation}
 P_{ab} = \sum_{k=0}^{N_y-1} \left[\sum_{j=0}^{N_x-1} p_{jk} e^{-2\pi \mathrm{i} aj/N_x}\right] e^{-2\pi \mathrm{i} bk/N_y}.
\end{equation}

\noindent We can view the expression inside the brackets as \(P_{ak}^x = \sum_{j=0}^{N_x-1} p_{jk} e^{-2 \pi \mathrm{i} aj/N_x}\), as the DFT along the \(x\)-axis for a fixed index \(k\). Once we compute the $N_y$ DFTs along the \(x\)-axis, we can write:

\begin{equation}
 P_{ab} = \sum_{k=0}^{N_y-1} P_{ak}^x e^{-2\pi \mathrm{i} kb/N_y}.
\end{equation}

\noindent Notice that this yields another expression representing \(N_x\) DFTs in one dimension, this time along the \(y\)-axis for index \(a\). Consequently, we can obtain the DFT of a function defined on a tw-dimensional domain by calculating the DFT in one dimension sequentially. If both \(N_x\) and \(N_y\) are powers of two, we can efficiently compute them using the FFT. 

Finally, notice that this procedure can be extended analogously to functions defined on an \(n\)-dimensional domain.

\section{Applications}
\label{sec:applications}

In this work we consider applications involving equations whose unknowns depend on two spatial coordinates, under this restriction we show how to solve elliptic, parabolic and hyperbolic equations.

\subsection{Example: Elliptic PDE}

The Poisson equation in \(n\) dimensions is expressed as:

\begin{equation}
    \nabla^2 u(\Vec{x}) = s(\Vec{x}), \qquad \Vec{x} \in D \subset \mathbb{R}^n,
    \label{eq:Poisson}
\end{equation}

\noindent where \(u\) is the potential function and \(s\) is a source term. This equation can be solved given appropriate boundary conditions on \(\partial D\). To analyze the Poisson equation in Fourier space, we apply the Fourier transform, denoted by \(\mathcal{F}\), to obtain:

\begin{equation}
    -\omega^2 \mathcal{F}\{u\} = \mathcal{F}\{s\},
    \label{eq:PoissonFourier}
\end{equation}

\noindent this relation is derived from the differentiation property of the Fourier transform, which states that:

\[
\mathcal{F} \left\{ \frac{\partial^k u(\Vec{x})}{\partial x_i^k} \right\} = (-\mathrm{i}\omega_i)^k \mathcal{F}\{u(\Vec{x})\}.
\]

\noindent where $x_i$ and $\omega_i$ are the $i-th$ coordinates of $\Vec{x}$ and $\Vec{\omega}$ in the coordinate and frequency spaces respectively. Foe the calculation of the FT of \(u\),  specific conditions are to be filfilled. One critical condition is the regularity of the source function \(s\). For the numerical methods employed to approximate the solution, it is necessary  \(\mathcal{F}\{u\}(\omega=0)\) to be finite. This requirement implies that the left-hand side of equation (\ref{eq:PoissonFourier}) must vanish, leading us to conclude that the right-hand side must also be zero.

To ensure this, we define a new function:

\[
g = s - \bar{s},
\]

\noindent where \(\bar{s}\) represents the average value of the source term \(s\) over the volume \(V\):

\begin{equation}
    \bar{s} = \frac{1}{V} \int s(\Vec{x}) d^n x,
\end{equation}

\noindent with \(V = \int_D d^n x\) being the volume occupied by the domain \(D\). Thus, we can rewrite Poisson equation as:

\begin{equation}
    \nabla^2 u = g = s - \bar{s}.
    \label{eq:Poisson2}
\end{equation}

\noindent The solution to the equation (\ref{eq:Poisson2}) can therefore be expressed as:

\begin{equation}
    u(\Vec{x}) = \mathcal{F}^{-1}\left\{ -\frac{\mathcal{F}\{g\}}{\omega^2} \right\},
    \label{eq:SolutionPoisson}
\end{equation}

\noindent However, this expression can only be solved exactly in specific cases. In practical applications, we replace the Fourier operator \(\mathcal{F}\) with the Fast Fourier Transform (FFT) in the discrete domain \(D_d\), sampled with resolutions \(\Delta x_i\).

As an example of the use of the FFT in this type of problem, consider the two-dimensional case $n=2$. More specifically let us consider the Poisson equation sourced by a Gaussian:

\begin{equation}
    s(x,y) = Ae^{-(x^2 + y^2)/\sigma^2},
\end{equation}

\noindent and the problem defined in the domain \(D = [-2,2]^2\), sampled with resolutions \(\Delta x = \Delta y = \frac{4}{N}\), where \(N = 64\) and \(N = 128\). Moreover consider the parameters \(A = 1\) and \(\sigma = 0.1\). 

\begin{figure}
    \centering
    \includegraphics[width=8cm]{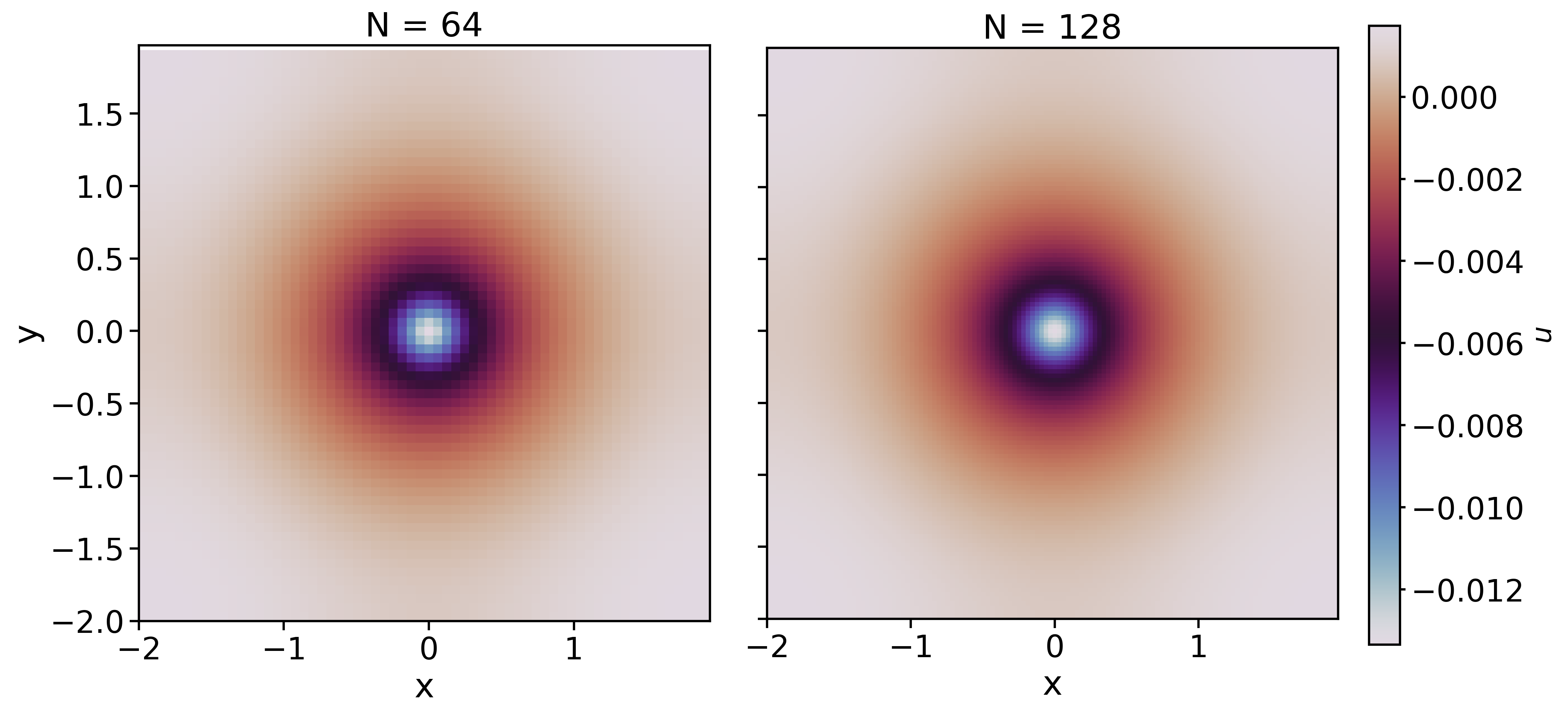}
    \caption{Numerical solutions of the two-dimensional Poisson equation with a Gaussian source. The left panel shows the solution for a grid resolution of \(N = 64\), while the right panel displays the solution for \(N = 128\). The results demonstrate how increasing the grid resolution leads to a more refined and accurate representation of the solution.}
    \label{fig:PoissonSolution}
\end{figure}

The numerical solution to this problem is shown in Figure \ref{fig:PoissonSolution}. The left panel presents the solution with \(N = 64\), while the right panel shows the solution with \(N = 128\), illustrating the consistency of the numerical solutions.

\subsection{Example: Parabolic equation}

The paradigm of a parabolic equation is the Diffusion Equation:

\begin{equation}
    \dfrac{\partial u}{\partial t} = \nabla^2 u + s, \quad \vec{x} \in D \subset \mathbb{R}^n,
    \label{ec:EcuacionDeCalor}
\end{equation}

\noindent which defines an initial value problem (IVP) provided the initial conditions for \(u\):

\begin{equation}
    u(0, \Vec{x}) = f(\Vec{x}).
    \label{ec:CalorCI}
\end{equation}

\noindent Here, \(s\) is a source term, typically interpreted as a reaction term. In Fourier space, this IVP becomes:

\begin{equation}
    \dfrac{d\hat{u}}{dt} = -\omega^2\hat{u} + \hat{s}, \quad \hat{u}(0) = \hat{f}(\Vec{\omega}),
    \label{eq:DiffusionIVPFS}
\end{equation}

\noindent where \(\hat{u} = \mathcal{F}(u)\) and \(\hat{s} = \mathcal{F}(s)\) are the FT of \(u\) and \(s\), respectively. Notice that this expression is a set of  Ordinary Differential Equations (ODEs) for $\hat{u}$ and for each value of $\omega$. This set of ODEs for different values of \(\omega\) can be solved numerically in general with simply an ODE integrator. Fortunately, for the special case where \(s=0\), which we also consider here, there is a closed-form solution:

\begin{equation}
    \hat{u} = \left\lbrace
    \begin{matrix}
    \hat{f} e^{-\omega^2 t} & \text{if } \omega>0, \\
    \hat{f} & \text{if } \omega=0,
    \end{matrix}
    \right.
    \label{ec:solCalorFourier}
\end{equation}

\noindent and therefore the solution to the original problem is \(u = \mathcal{F}^{-1}(\hat{u})\). In most cases, this cannot be solved exactly, but we can observe the asymptotic behavior: as \(t \to \infty\), the unique value that does not vanish is \(\hat{f}(0)\), leading to the conclusion that the asymptotic solution of the Diffusion Equation is \(u \to \bar{f}\). For example, in the context where \(u\) represents temperature, it redistributes throughout the entire domain, eventually reaching thermodynamic equilibrium, consequently the quantity \(\bar{u}\) must be considered using the FFT method.

\begin{figure}
    \centering
    \includegraphics[width=8cm]{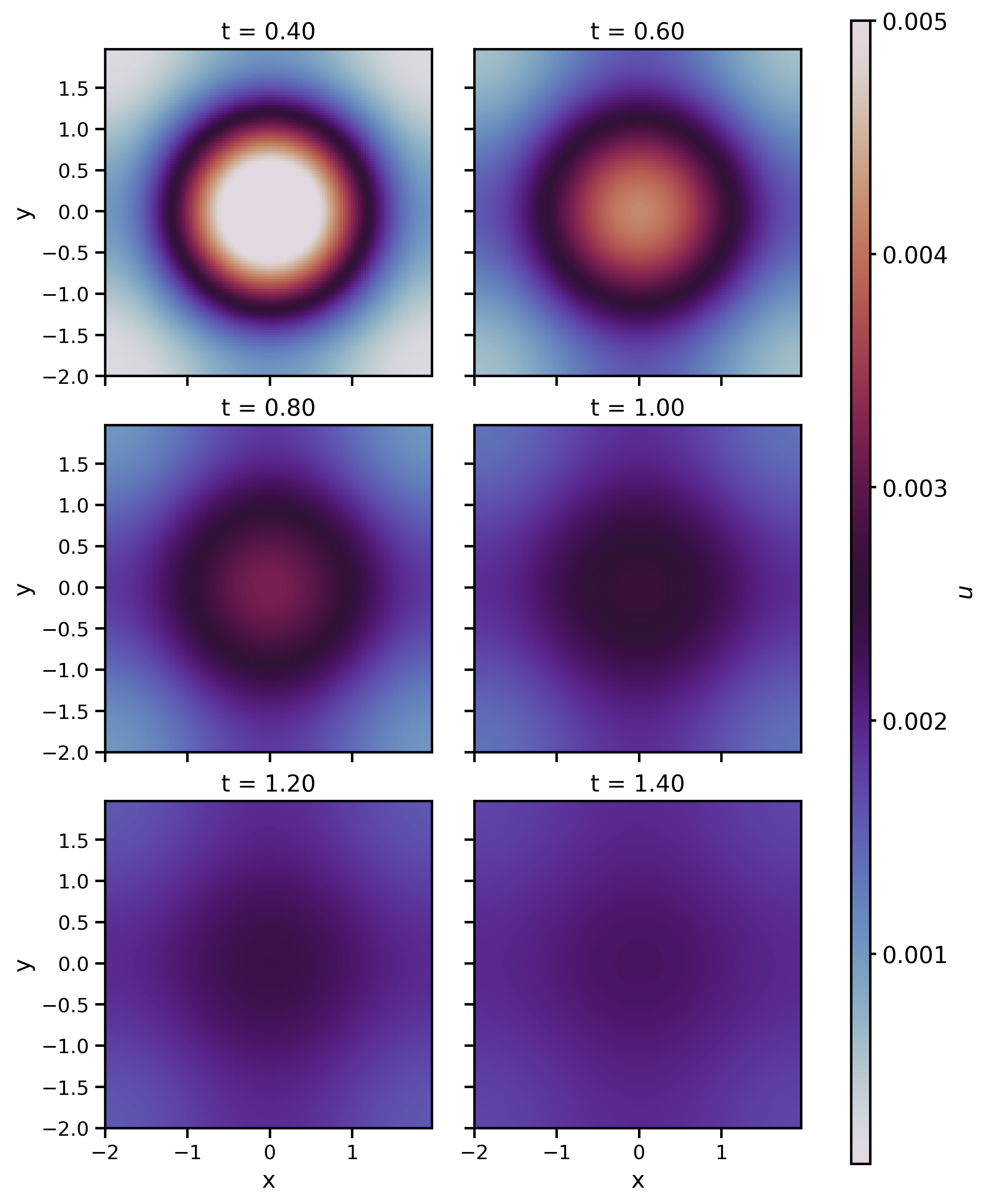}
    \caption{Snapshots of the numerical solution of the two-dimensional Diffusion Equation with a Gaussian initial condition without source, shown at different times. The diffusion process leads to a progressive dissipation of the initial profile.}
    \label{fig:diffusionsnaps}
\end{figure}

As a particular example in two dimensions, we solve the Diffusion Equation in the spatial domain \(D = [-1,1]^2\), discretized with \(\Delta x = \Delta y = \frac{2}{N}\), where \(N = 128\) points along each direction. The time domain is set as \(t \in [0,10]\). We solve this IVP, without source $s=0$, using the following initial Gaussian profile:

\begin{equation}
    f(x,y) = A e^{-(x^2 + y^2)/\sigma^2},
    \label{eq:f_1}
\end{equation}

\noindent where \(A = 1\) and \(\sigma = 0.1\). Figure \ref{fig:diffusionsnaps} shows snapshots of the solution at various time steps. As expected, the initial Gaussian profile dissipates over time, gradually approaching a uniform value, consistent with the asymptotic behavior of the diffusion process.

\begin{figure}
    \centering
    \includegraphics[width=8cm]{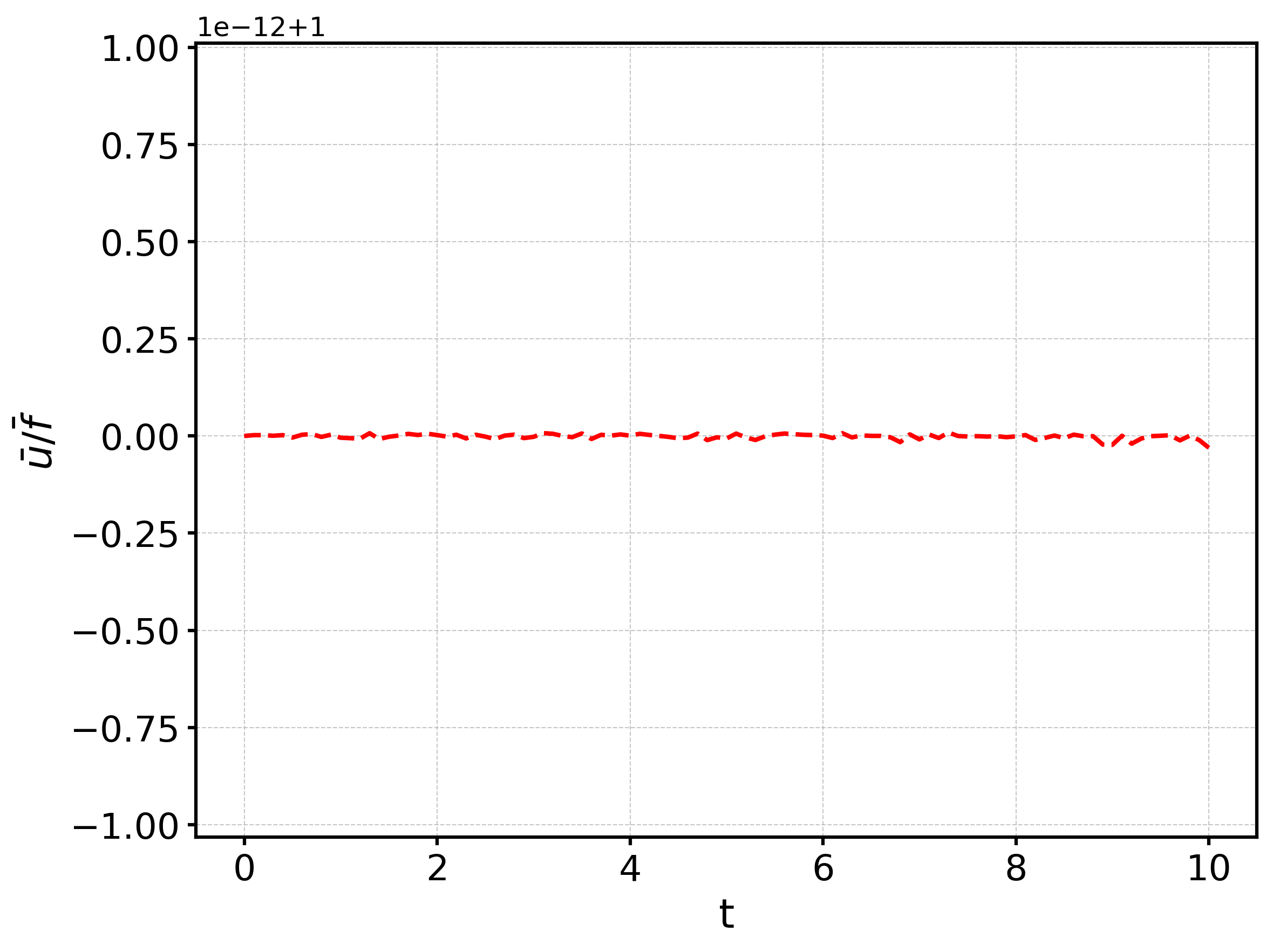}
    \caption{Normalized average value of the function \(u\) over time, relative to the initial average \(\bar{f}\). The plot highlights the conservation of \(\bar{u}\) throughout the diffusion process. Notice that the scale in the vertical axis is $10^{-12}+1$, thus the plot indicates that the ratio is nearly a perfect $\bar{u}/\bar{f}\sim 1$.}
    \label{fig:diffusionerror}
\end{figure}

In Figure \ref{fig:diffusionerror}, we present the time evolution of the average value \(\bar{u}\), normalized with respect to \(\bar{f}\), the average of the initial condition. The result demonstrates that the FFT method conserves this quantity throughout the evolution.

\subsection{Example: Hyperpolic PDE}

The study of the wave equation is a classic in Physics and Mathematics, because it describes a wide variety of phenomena, from the propagation of light and sound until Einstein's equations for the evolution of the spacetime geometry \cite{guzman2010}. The non-homogeneous wave equation is a hyperbolic partial differential equation, and in $n$ dimensions, it is written as

\begin{equation}
\dfrac{\partial^2u}{\partial t^2} - c^2\nabla^2u = s, \quad \vec{x} \in D \subset  \mathbb{R}^n,
    \label{eq:wave_equation}
\end{equation}

\noindent where $D$ is the domain of the wave function $u$,  $s$ is a source which in general can depend on position and time and $c$ is the velocity of propagation. In what follows we use $c=1$, knowing that the solutions for another arbitrary velocity can be found by rescaling time as $t\rightarrow ct$.

Equation (\ref{eq:wave_equation}) can be viewed as an Initial Value Problem (IVP) provided initial conditions for $u$ and $\partial_t u$ at initial time $t=0$ and boundary conditions at $\partial D$:

\begin{eqnarray}
u(0,\Vec{x}) &=& f(\Vec{x}) \nonumber\\
\dfrac{\partial u}{\partial t} (0,\Vec{x}) &=& g(\Vec{x}).\label{ec:OndaCI}
\end{eqnarray}

\noindent As mentioned above, the FT method can be easily applied in periodic domains, thus boundary conditions for the problem (\ref{eq:wave_equation})-(\ref{ec:OndaCI}) will be assumed to be periodic. These ingredients define a well posed IVP. 

For the construction of the solution we remind that $\hat{u} = \mathcal{F}(u)$. Then, equation (\ref{eq:wave_equation}) and initial conditions (\ref{ec:OndaCI}) in the Fourier space read:

\begin{equation}
\dfrac{d^2\hat{u}}{dt^2}+\omega^2\hat{u} = \hat{s}, \text{ I.C. }\left\lbrace
 \begin{matrix}
     \hat{u}(0) = \hat{f}(\Vec{\omega})   \\ \\
\dfrac{d\hat{u}}{dt}(0) = \hat{g}(\Vec{\omega}).
 \end{matrix}
\right.
\label{eq:WaveEquationFourierSpace}
\end{equation}

\noindent Notice that the wave equation in Fourier space is a family of ODEs for forced harmonic oscillators driven by an external force $\hat{s}$, for each value of the frequency $\omega$. In order to use standard ODE integrators we write the oscillator equations as first order systems of equations:

\begin{eqnarray}
    \dfrac{d\hat{u}}{dt} &=& v,\nonumber\\
    \dfrac{dv}{dt} &=& -\omega^2\hat{u}+\hat{s},
\label{eq:SystemOfODEsForWaveEquation}
\end{eqnarray}

\noindent In general these equations can be integrated numerically using for example any flavor of Runge-Kutta methods. In the case when $s$ is time independent, the system reduces to that of a forced oscillator and exact solutions are found:

\begin{itemize} 
 \item[(i)] For $\omega=0$

\begin{equation}
\hat{u}(t) = \dfrac{1}{2}\hat{s}t^2+\hat{g}t+\hat{f},
    \label{eq:OmegaZero}
\end{equation}

\item[(ii)] For $\omega\neq0$

\begin{equation}
\hat{u}(t) = (\hat{f}-\dfrac{\hat{s}}{\omega^2})\cos{(\omega t)}+ \dfrac{\hat{g}}{\omega}\sin(\omega t).
    \label{eq:OmegaNonZero}
\end{equation}
\end{itemize}

\noindent Then, the only step remaining to construct the solution to the original wave equation is to calculate $u = \mathcal{F}^{-1}(\hat{u})$.

In summary, the wave equation in a periodic domain is solved using the FFT with the following recipe:

\begin{itemize}
    \item[-] Define the equation and initial conditions.
    \item[-] Write the equation and initial conditions in Fourier space.
    \item[-] Write the resulting equation as a set of first-order ODEs in time.
    \item[-] Solve the system using numerical integration. If $s$ is time-independent, use expressions (\ref{eq:OmegaZero}) and (\ref{eq:OmegaNonZero}) to find the solution $\hat{u}$.
    \item[-] Recover the solution to the original equation using the IFFT.
\end{itemize}

Now, we exemplify the method for the case of two spatial dimensions with two cases: the first one involves a source term \(s = 0\), and the second one includes a time-dependent source.

\subsubsection{Example with $s=0$.}

\begin{figure}
\centering
\includegraphics[width=8cm]{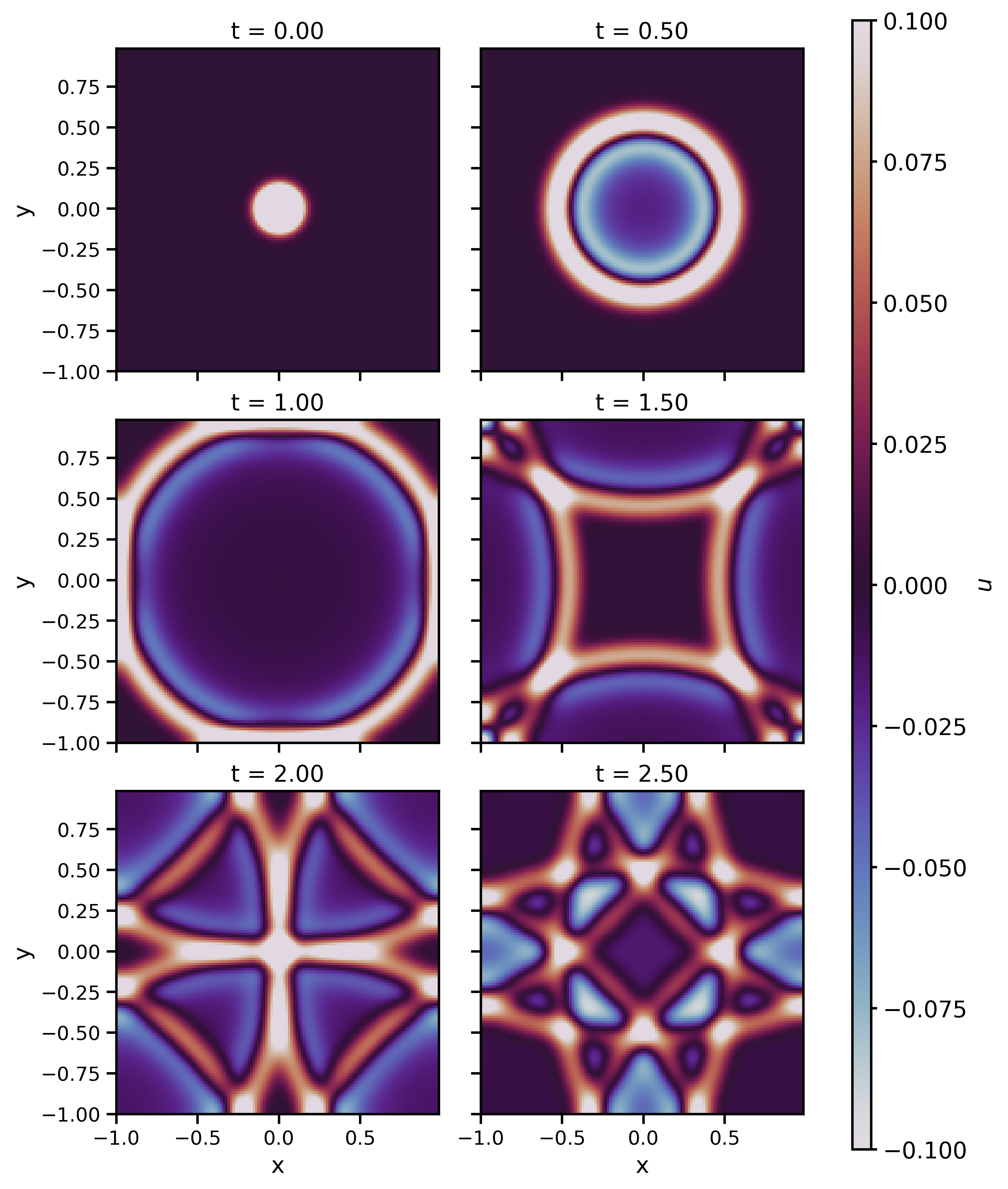}
\caption{Snapshots of the solution for the case $s=0$, using $N = 128$ at times $t=0,~0.5,~1,~1.5,~2.0$ and 2.5.}
\label{fig:solution64}
\end{figure}

The problem is to solve the wave equation for an initial time-symmetric spherically symmetric Gaussian pulse in the domain $D=[-1,1]^2$. The set up of the problem is then with $s=0$ and initial conditions

\begin{eqnarray}
f(x,y) &=& Ae^{-(x^2+y^2)/ \sigma^2}\nonumber\\
g(x,y) &=& 0.\nonumber
\end{eqnarray}

\noindent in the discrete domain that uses $N=128$, to be integrated in the time domain $t\in[0,2]$. In Figure \ref{fig:solution64} we illustrate the numerical solution at various times using the equations (\ref{eq:OmegaZero}) and (\ref{eq:OmegaNonZero}). Notice that the initial Gaussian spherical pulse evolves as a spherical wave whose amplitude decreases as $1/r$ prior to reaching the boundary, when the wave-front arrives to the boundary of the domain, it reenters through opposite side of the domain and starts to produce the expected interference patterns.

\begin{figure}
    \centering
    \includegraphics[width=8cm]{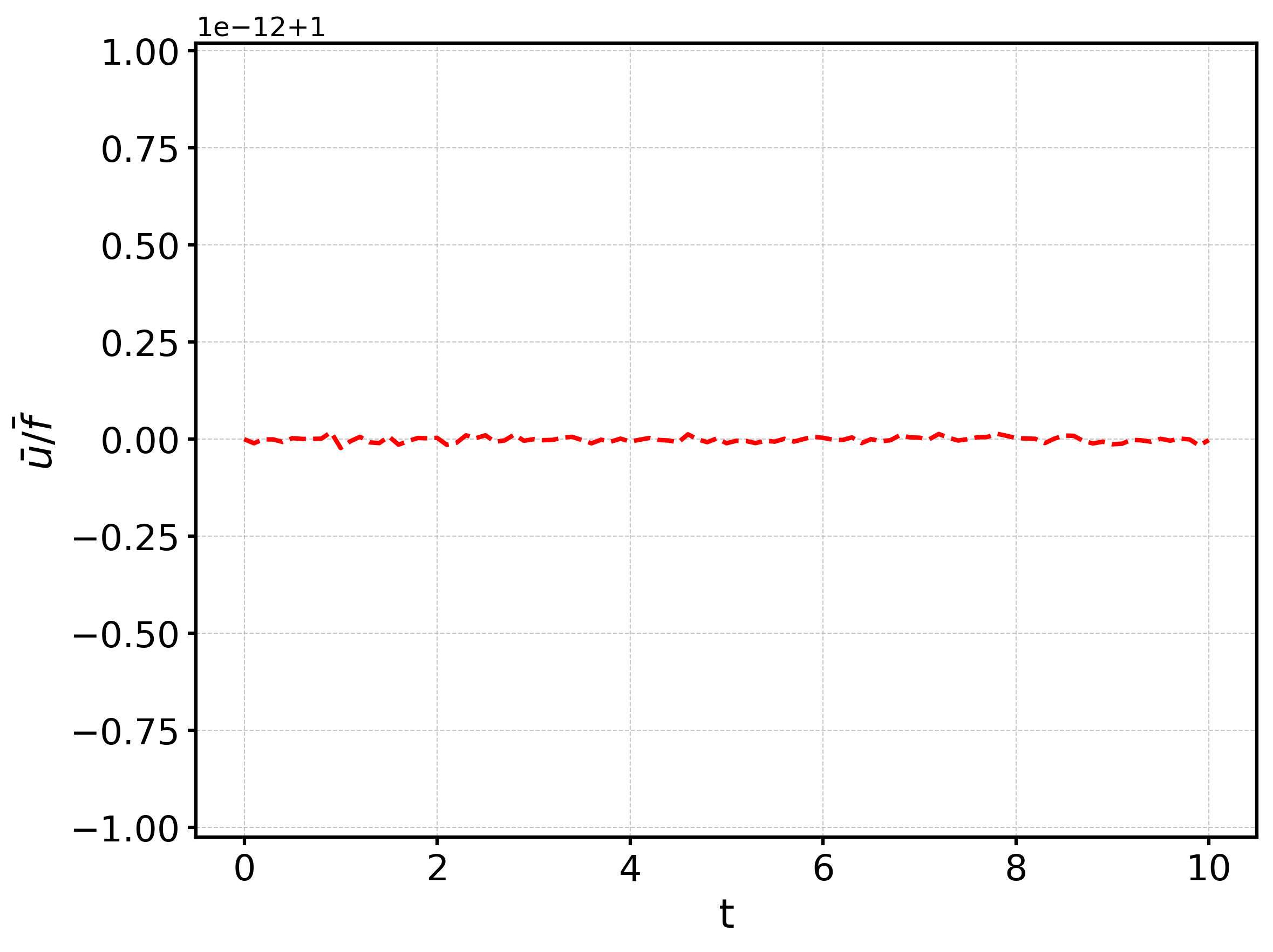}
    \caption{Normalized average value of the function \(u\) over time, relative to the initial average \(\bar{f}\). The plot highlights the conservation of \(\bar{u}\) throughout the evolution. Notice again, that the scale in the vertical axis is $10^{-12}+1$, thus the plot indicates that the ratio is nearly a perfect $\bar{u}/\bar{f}\sim 1$.}
    \label{fig:wave1error}
\end{figure}

In Figure \ref{fig:wave1error}, we present the time evolution of the average value \(\bar{u}\), normalized with respect to \(\bar{f}\), the average of the initial condition. The results demonstrate that the FFT method conserves this quantity  throughout the simulation.

\subsubsection{Example with Source}

\begin{figure}[h]
    \centering
    \includegraphics[width=8cm]{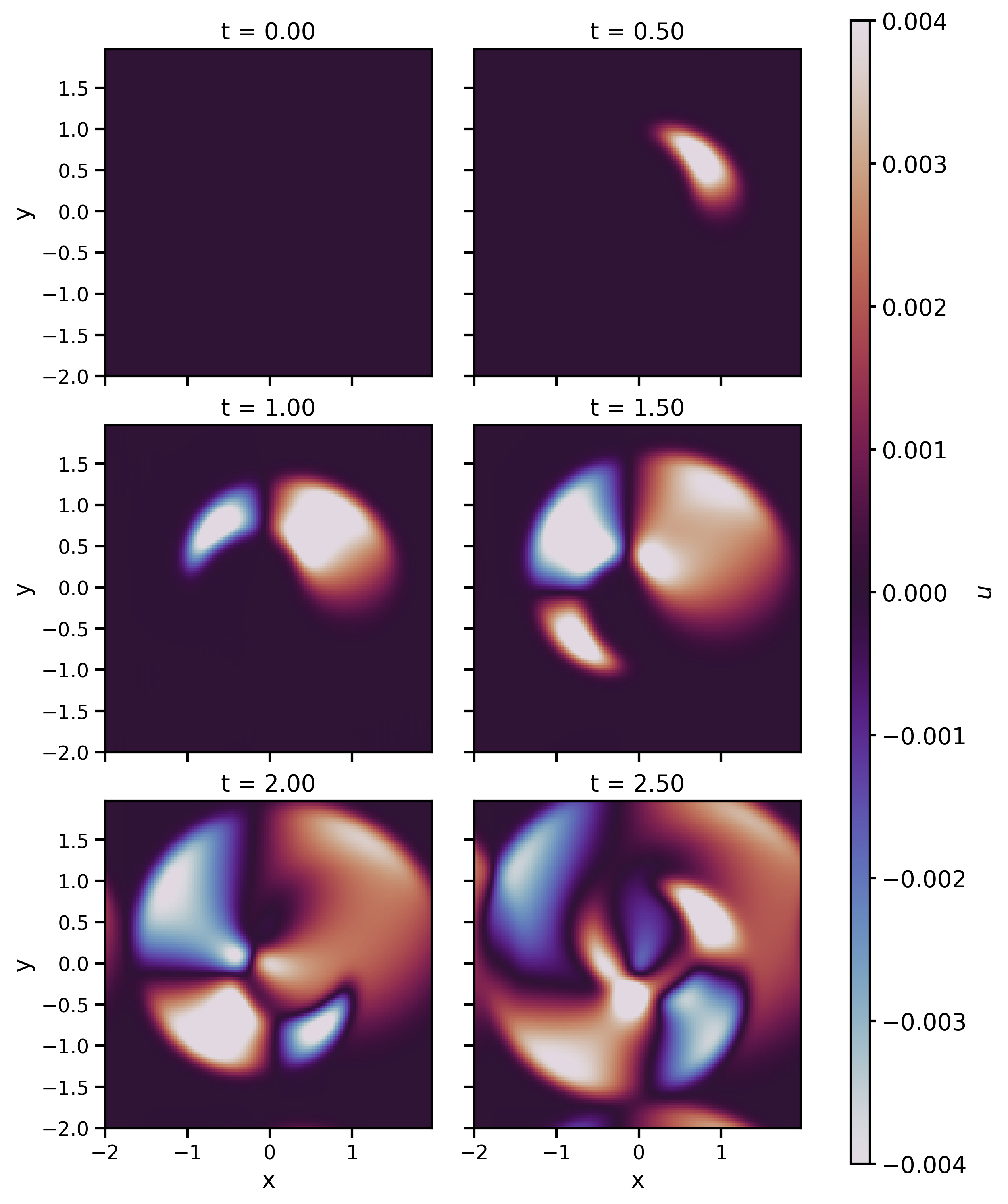}
    \caption{Snapshots of the solution to the wave equation with an orbiting source, computed using \(N = 128\) at times \(t = 0,~0.5,~1,~1.5,~2.0,\) and \(2.5\). The solution displays the expected evolution over time.}
    \label{fig:wesources}
\end{figure}

Another typical solution of the wave equation involves a source that orbits the origin on the \(xy\)-plane, leaving behind a strip. The source is given by:

\begin{equation}
    s(t, x, y) = A e^{-[(x - x_0)^2 + (y - y_0)^2]/\sigma_s^2} \cos(\gamma t),
\end{equation}

\noindent where the center of the source \((x_0, y_0)\) orbits in a circle of radius \(r_s\), described by \(x_0 = r_s \cos(\Omega t)\) and \(y_0 = r_s \sin(\Omega t)\). In this example, we use the parameters \(A = 1\), \(\sigma_s = 0.1\), angular frequency \(\Omega = 5\), and \(\gamma = 10\), within the domain \(D = [-2, 2]^2\) discretized with spatial resolution \(h = \Delta x = \Delta y = 4/N\), where \(N = 64\), 128, and 256. The temporal resolution is set to \(\Delta t = 0.25 h\). Various snapshots of the solution, obtained using the Runge-Kutta 4 (RK4) scheme to integrate Eq. (\ref{eq:SystemOfODEsForWaveEquation}), are shown in Figure \ref{fig:wesources}.

Next, we perform an aself-convergence test using three numerical solutions. Let \(u_1\), \(u_2\), and \(u_3\) represent the solutions obtained with resolutions \(h\), \(h/2\), and \(h/4\), respectively; notice that since the spatial and temporal resolutions are proportional to \(h\), time resolution is being also doubled. The exact solution \(u_0\) can be approximated as:

\[
    u_0 = u_l + \epsilon \left(\frac{h}{2^{l-1}}\right)^4,
\]

\noindent with \(l = 1, 2, 3\), where the exponent 4 corresponds to the fourth-order error of the RK4 method. Thus, if the numerical solutions self-converge, the following relation between the numerical solutions should hold:

\begin{equation}
    \frac{u_2 - u_1}{u_3 - u_2} \sim 2^4.
    \label{eq:SCF}
\end{equation}

\begin{figure}[h]
    \centering
    \includegraphics[width=8cm]{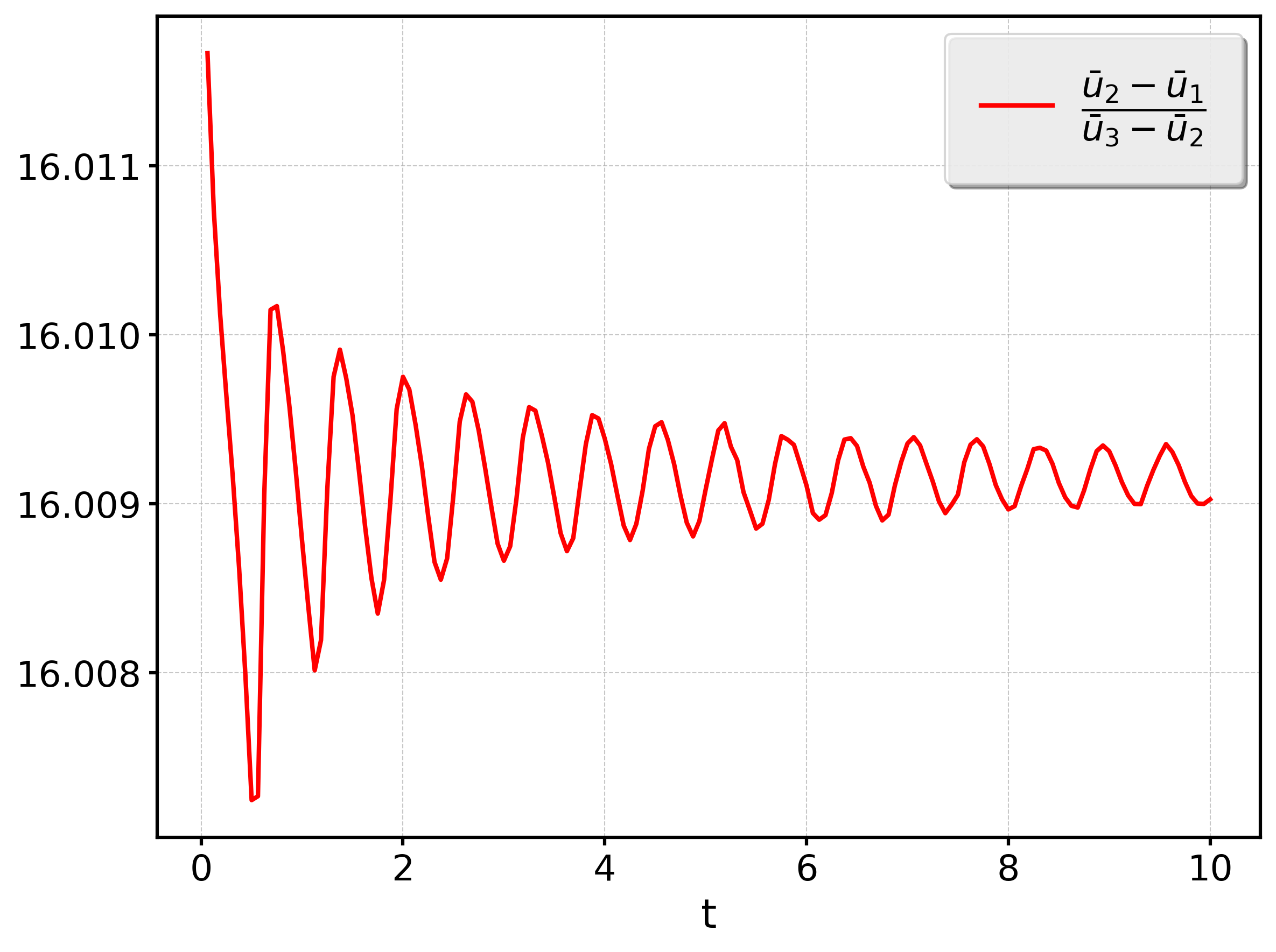}
    \caption{Self-convergence test for the numerical solutions of the wave equation with an orbiting source, calculated using \(N = 64\), \(128\), and \(256\) points. The results show a convergence factor of approximately 16=$2^4$, as expected for the RK4 method.}
    \label{fig:ACT}
\end{figure}

\noindent Although this is a local condition, we can verify it by examining global quantities, such as the average values \(\bar{u}_l\). Figure \ref{fig:ACT} confirms that the relation above holds beautifully. This indicates that the solutions are self-convergent, in one word, numerically correct.

\section{Final comments}
\label{sec:conclusions}

We have presented the basic elements needed for the calculation of the DFT and the FFT in one and more dimensions, along with the restrictions on the calculations of the this transform. We show explicitly with an example the scaling of the complexity of the two algorithms, illustrating the typical ${\cal O} (N\log N)$ number of operations of the FFT, which is smaller than the ${\cal O} (N^2)$ complexity of the DFT when the number of points of the discrete domain $N$ is large.

We then applied the calculation of the FFT to the solution of Partial Differential Equations, in examples involving two spatial dimensions. Specifically, Poisson equation to illustrate the solution of an elliptic equation, the diffusion equation as example of a parabolic equation and the wave equation as the paradigm of hyperbolic equations. We have shown step by step how to construct reliable numerical solutions of these examples and expect this paper sets a starting point for students to tackle state of the art problems with this method.

As motivation to lear this method, we can mention that the method of this paper serves to solve systems of coupled equations, for example the simultaneous Schr\"odinger-Poisson system that helps modeling bosonic dark matter \cite{Mocz:2017wlg,CAFE2022}, and is also the workhorse method to simulate stationary and dynamic scenarios of laboratory Bose-Einstein Condensates \cite{BECs}. Going even further, the quantum version of the FT is currently been implemented to solve PDEs in quantum computers \cite{GuzmanGuzman2024}, including an interesting variety of physical scenarios like plasma physics and fluid dynamics \cite{AppPlasmaPhysics,AppHeatEquation,AppNavierStokes}, as well as financial applications \cite{BlackScholes}.


\section*{Acknowledgments}
Iv\'an \'Alvarez receives support within the CONACyT graduate scholarship program under the CVU 967478. This research is supported by grants CIC-UMSNH-4.9 and Laboratorio Nacional de C\'omputo de Alto Desempe\~no Grant No. 1-2024.


\bibliography{FFT}

\end{document}